\documentclass{article}
\usepackage[]{epsf,epsfig,amsmath,amssymb,amsfonts,latexsym}
\title{ Tail Invariant Measures of the Dyck Shift}
\author{Tom Meyerovitch}

\newtheorem{theorem}{Theorem}[section]
\newtheorem{lemma}{Lemma}[section]
\newtheorem{proposition}[theorem]{Proposition}

\newtheorem{remark}{Remark}[section]

\newenvironment{proof}{{\bf Proof:} \rm}{\hfill $\Box$ \medskip\\}

\begin{document}
\maketitle
\begin{abstract}
We show that the one-sided Dyck shift has a unique tail invariant
topologically $\sigma$-finite measure (up to scaling). This
invariant measure of the one sided Dyck turns out to be a
shift-invariant probability. Furthermore, it is one of the two
ergodic probabilities obtaining maximal entropy. For the two sided
Dyck shift we show that there are exactly three ergodic
double-tail invariant probabilities. We show that the two sided
Dyck has a double-tail invariant probability, which is also shift
invariant, with entropy strictly less than the topological
entropy.\footnote{This article is a part of the author's M.Sc.
thesis, written under the supervision of J. Aaronson, Tel-Aviv
University.}
\end{abstract}

\section{Introduction}
The study of tail invariant probabilities for subshifts has so far
focused mostly on sofic systems. There are known results for the
case of the one sided tail of (mixing) SFT's \cite{BM}. Also, for
the case of the $\beta$-shift it is known that there exists a
unique tail-invariant measure \cite{ANS}. In all of these examples
the tail-invariant measure is also equivalent to a unique shift
invariant measure of maximal entropy. Invariant measures for the
double-tail (and some sub-relations of the double-tail) of SFT's
have also been characterized \cite{PS97}.\\
 Let $\Sigma$ be a finite alphabet. For a subshift $X \subset \Sigma^\mathbb{Z}$, we
define the double-tail relation, or \emph{homoclinic} \cite{PS97}
relation of $X$ to be:
\[\mathcal{T}_2(X) := \{(x,x') \in X \times X\:\ \exists n \geq 0 \ \forall |k|>n \ x_k=x'_k \}\]
A $\mathcal{T}_2(X)$-holonomy is an injective Borel function $g:A
\mapsto g(A)$, with $A$ a Borel set and $(x,g(x)) \in
\mathcal{T}_2(X)$ for every $x \in A$. We say that $\mu \in
\mathcal{M}(X)$ is a double-tail invariant if $\mu(A)=\mu(g(A))$
for every $\mathcal{T}_2(X)$-holonomy $g$.\\
In this paper we identify the tail invariant probability measures
for the Dyck Shift. This subshift was used in \cite{WK74} as a
counter-example for a conjecture of B. Weiss, showing there are
exactly two measures of maximal entropy for this subshift, both of
which are Bernoulli. We show that for the one-sided Dyck shift one
of these measures is the unique tail-invariant probability
(section \ref{dyck1}). We also characterize the double-tail
invariant probabilities for the Dyck shifts (section \ref{dyck2}).
In addition to its two equilibrium measures, the two sided Dyck
shift has another double-tail invariant probability -- shift
invariant, non-equilibrium. These are the only three double-tail
invariant, ergodic probabilities on the two sided Dyck shift. A
different but perhaps related study of the Dyck shift was carried
out by Hamachi and Inoue \cite{HI05}.

\section{Definition of The Dyck System}

Let us explicitly describe the Dyck language and it's
cover (Fischer automaton).\\
 Let $m \geq 1$ and
\(\Sigma=\{\alpha_j : 1 \leq j \leq m \}\cup \{\beta_j: 1 \leq j
\leq m\}\), \(\Gamma=\{\alpha_j: 1\leq j \leq m \}^*\), and with
$\Lambda$ the empty word, $\varphi(a,\alpha_j)=a\alpha_j$,$a\in
\Gamma$,
\[\varphi(a,\beta_j)=\left\{
\begin{array}{ll}
\beta_j & \mbox{if $a=\Lambda$,or $a=(a_j)$} \\
\Lambda & \mbox{if $a\in\{a_j: 1 \leq j \leq m\}^k$,$k\in
\mathbb{N}$, $a_k \neq \alpha_j$} \\
(a_i)_{i=1}^{k-1} & \mbox{if $a\in \{a_j: 1 \leq j \leq m\}^k$,
$k>1$,$a_k=\alpha_j$}
\end{array}
\right.
\]
 Another way to describe the Dyck-Shift is in terms of it's syntactic monoid:\\

Let $M$ be the monoid generated by $\Sigma$, with the following
relations:\\
\begin{enumerate}
\label{def_m}
  \item \(\alpha_j \cdot \beta_j \equiv \Lambda \equiv 1 (mod M) , j=1, \ldots ,m\)
  \item \(\alpha_i \cdot \beta_j \equiv 0 (mod M) ,i \neq j \)
\end{enumerate}
The \textsl{m-Dyck language} is
\[ L=\{l \in \Sigma^* : l \neq 0 (mod M) \} \]
and the corresponding \text{(two sided) m-Dyck subshift} is
\[ X=\{ x \in \Sigma^\mathbb{Z} :\; (x_i)_{i=r}^{l} \in L \mbox{ for
all } -\infty  <r \leq l < +\infty \} \] and we will also refer to
the \textsl{one sided m-Dyck subshift}:
\[ Y=\{ y \in \Sigma^\mathbb{N} :\; (y_i)_{i=r}^{l} \in L \mbox{ for
all } 0  \leq r \leq l < +\infty \} \]

These are indeed subshifts, since we only pose restrictions on
finite blocks. Conversely, we will later note by $L(X)=L(Y)=L$ the
language consisting of words which are admissible in $X$. Also,
let:
\[\mathcal{L}_n= \mathcal{L}(Y,n):=L(Y) \cap \Sigma^n\]

Note that when $m=1$, $X$ is simply the full 2-Shift, and so we
will only be interested in the case where $m\geq 2$.\\

For $w=(w_0,\ldots,w_{n-1}) \in \mathcal{L}(X,n)$ define
\[H(w)=\sum_{i=0}^{n-1}\sum_{j=1}^{m}(\delta_{\alpha_j,w_i}-
\delta_{\beta_j,w_i})\] and $H(\Lambda)=0$.

 For $x \in X$, let
\begin{equation}
\label{H_def}
 H_i(x)= \left\{ \begin{array}{ll}
  \sum_{j=0}^{i-1}
\sum_{l=1}^{m}(\delta_{\alpha_l,x_j}-\delta_{\beta_l,x_j}) &
\mbox{if
$i>0$} \\
\sum_{j=i}^{-1}\sum_{l=1}^{m}(\delta_{\beta_l,x_j}-\delta_{\alpha_l,x_j})
& \mbox{if $i<0$} \\
0 & \mbox{if $i=0$} \\
\end{array}
\right.
\end{equation}

We shall use the same notation for the one-sided subshift. For $y
\in Y$, let
\[ H_i(y)= \left\{ \begin{array}{ll}
\sum_{j=0}^{i-1}\sum_{l=1}^{m}(\delta_{\alpha_l,y_j}-\delta_{\beta_l,y_j})
& \mbox{if $i>0$} \\
0 & \mbox{if $i=0$} \\
\end{array}
\right.
\]
 where it is clear from the context whether we are
 refereing to the one sided or two sided subshift.
If $w \equiv 1$ (mod $M$) we say that $w$ is a \emph{balanced
word}.\\
A word $w $ is a \emph{Dyck word} if  it is a minimal balanced
word. This means $w=\alpha_i \tilde{w} \beta_i$ for some balanced
word $\tilde{w}$ and $1 \le i \le m$.\\

For $w \in \mathcal{L}_n$ define
$$m(w)=\min\{H(u):\; u \mbox{ is a prefix of }w\}$$
$$\hat{\alpha}(w):=H(w)-m(w)$$
$$\hat{\beta}(w):=-m(w)$$
where in the definition of $m(w)$ it is understood that the empty
word is a prefix of any word, so that $m(w) \le 0$. $\hat{\alpha}$
is the number of unmatched $\alpha$'s in $w$, and $\hat{\beta}$ is
the number of unmatched $\beta$'s in $w$. We say that $w$ has an
unmatched $\alpha$ at location $t$ if $w_t=\alpha_i$, and
$\hat{\alpha}(w_{0}^{t-1})< \hat{\alpha}(w_{0}^{t})$. We define
''unmatched $\beta$'' respectively using $\hat{\beta}$. We say
that $x \in X$, has an unmatched $\alpha$ ($\beta$) at location
$t$ if $x_t=\alpha_i$ ($x_t=\beta_i$) which is unmatched in any
finite word $x_{[a,b]}$ with $t\in [a,b]$.

 \subsection{Classification of the Dyck System}
 Before stating and proving the result regarding invariant
 measures for the Dyck system, we characterize this subshift in
 terms of the classes of subshifts introduced in \cite{BH86},\cite{FI92}, and
 \cite{B05}.
  The purpose of this subsection is to put in broader context the Dyck shift and the results in the following
  sections. Detailed discussions of these classes of subshifts can be found in the references above.   \\
 By defining the Dyck language as the language recognized by a
 Fischer automaton, we showed that the Dyck system is a coded
 system (as in \cite{BH86}). We claim that the Dyck system is half-synchronized, yet
 not synchronized (as in \cite{FI92}):
 \begin{proposition}\label{dyck_half_synch}
 Every word $w$ in the Dyck langauge is half synchronizing
 \end{proposition}
 \begin{proof}
 Suppose $w=w_0,\ldots,w_{n-1}$
 Let $(u_k)_{k \in \mathbb{N}}$ be an enumeration of the Dyck words. We define a left infinite sequence
 $x \in X$ as the word $w$ (ending in coordinate $0$), preceded by a concatenation of the words $(u_k)_{k \in
 \mathbb{N}}$, and followed by an infinite sequence of $\alpha_j$'s. $x$ is a left-transitive point.
 $\omega_{+}(x(-\infty,0])=\omega{+}{w}$, since every unmatched $\alpha_j$
 in $x(-\infty,0)]$ must be in $w$.
 \end{proof}

 \begin{proposition}
 The $m$-Dyck system is not a synchronized system, for $m >1$.
 \end{proposition}
 \begin{proof}
 Let $w \in L(X)$. There exist $l,r$ such that $lwr \equiv 1$ (mod
 $M$). Thus, for $i \neq j$, $\alpha_ilwr\beta_j \not\in L(X)$,
 but $\alpha_ilw \in L(X)$ and $wr\beta_j \in L(X)$. This show
 that $w$ is not a synchronizing word.
 \end{proof}

In \cite{B05} Buzzi, defined and studied a class of subshifts
called \emph{subshifts of quasi-finite type}. We state without
proof the following:
\begin{proposition}
For $m>1$, the $m$-Dyck system is not weak quasi-finite type.
\end{proposition}

\subsection{Maximal Measures for the Dyck Shift}
 In \cite{WK74} Krieger introduced
the following decomposition of $X$ into shift invariant subsets:
\[A_{+}=\{y\in X :
\lim_{i\to\infty}H_{i}(y)=-\lim_{i\to -\infty}H_{i}(y)=\infty\}\]
\[A_{-}=\{y\in X :
-\lim_{i\to\infty}H_{i}(y)=\lim_{i\to -\infty}H_{i}(y)=\infty\}\]
\[A_{0}=\bigcap_{i=-\infty}^{\infty}(\bigcup_{l=1}^{\infty}\{y\in X :
H_{i}(y)=H_{i+l}(y)\} \cap\bigcup_{l=1}^{\infty}\{y\in X :
H_{i}(y)=H_{i-l}(y)\})\]
 Since the complement of these sets, \( X
\backslash (A_+ \cup A_- \cup A_0) \) is a countable union of
wandering sets, every ergodic shift-invariant probability measure
assigns probability one to exactly one of these sets.\\
Let further
\[B_+=\bigcap_{i=-\infty}^{\infty}(\bigcup_{l=1}^m(\{x \in X: x_i=\alpha_l\} \cup
\bigcup_{k=1}^{\infty}\{x \in X:
x_i=\beta_l,H_{i-k}(x)=H_i(x)\}))\]
\[B_-=\bigcap_{i=-\infty}^{\infty}(\bigcup_{l=1}^m(\{x \in X: x_i=\beta_l\} \cup
\bigcup_{k=1}^{\infty}\{x \in X:
x_i=\alpha_l,H_{i+k}(x)=H_i(x)\}))\] and observe that $A_+ \cup
A_0 \subset B_+$, $A_- \cup A_0 \subset B_-$. Let $\Omega=
\{\alpha_1 \ldots \alpha_m,\beta\}^{\mathbb{Z}}$. Define \(
\widehat{H}_{0}(x)=0 \), \( \widehat{H}_i(x) =
\sum_{j=0}^{i-1}(\sum_{k=1}^{m}\delta_{x_j,\alpha_k}-\delta_{x_j,\beta})
\), $x\in \Omega$. Denote
\[\widehat{B}_+=\bigcap_{i=-\infty}^{\infty}(\bigcup_{l=1}^m\{\omega \in \Omega: \omega_i=\alpha_l\} \cup
\bigcup_{k=1}^{\infty}\{\omega \in \Omega:
\omega_i=\beta,\widehat{H}_{i-k}(\omega)=\widehat{H}_i(\omega)\})\]
\[\widehat{A}_{+}=\{\omega \in \Omega :
\lim_{i\to\infty}\widehat{H}_{i}(\omega)=-\lim_{i\to
-\infty}\widehat{H}_{i}(\omega)=\infty\}\]
\[\widehat{A}_{0}=\bigcap_{i=-\infty}^{\infty}(\bigcup_{l=1}^{\infty}\{\omega \in \Omega :
\widehat{H}_{i}(\omega)=\widehat{H}_{i+l}(\omega)\}
\cap\bigcup_{l=1}^{\infty}\{\omega\in \Omega :
\widehat{H}_{i}(\omega)=\widehat{H}_{i-l}(\omega)\})\]

 \[ (g_+(y))_{i}= \left\{
\begin{array}{ll}
\alpha_j & y_{i}= \alpha_j  \\
\beta & y_{i} \in \{ \beta_{1},\ldots,\beta_{m}\}
\end{array}
\right.
\]
$g_+$ is a one-to-one Borel mapping from $B_+$ onto
$\widehat{B}_+$, commuting with the shift. This shows that every
shift invariant probability measure $\mu$ on $X$ such that
$\mu(B_+)=1$ can be transported to a shift invariant probability
on $\Omega$ with equal entropy. By the intrinsic ergodicity of the
full-shift, there is a unique measure $\mu_1$ of maximal entropy
on $X$ such that $\mu_1(B_+)=1$. This measure is supported by $A+
\subset B_+$. By similar arguments, there is a unique measure
$\mu_2$ of maximal entropy on $X$ such that
$\mu_2(B_-)=1$, and in fact $\mu_2(A_-)=1$.\\
\begin{remark}
\[\sup_{\mu \in \mathcal{P}(A_0,T)}\{h(A_0,T,\mu)\} = \log(2)+\frac{1}{2}\log{m}\]
\end{remark}
\begin{proof}
 Since $A_0 \subset B_+$, any shift invariant probability $\mu_0$ on $X$ supported by
$A_0$ can also be transported to a probability $\widehat{\mu}_0$
on $\Omega$ via $g_+$. By the ergodicity,
$\widehat{\mu}_0([\beta])=\frac{1}{2}(1-\lim_{n \rightarrow
\infty}\frac{H_n(x)}{n})=\frac{1}{2}$. Thus,
$h(X,T,\mu_0)=h(\Omega,T,\widehat{\mu}_0) \leq
\log(2)+\frac{1}{2}\log(m)$, and equality can be obtained by
taking
$\widehat{\mu}_0=\prod_{i=-\infty}^{+\infty}(\frac{1}{2m},\ldots,\frac{1}{2m},\frac{1}{2})$.
\end{proof}

\section{Tail Invariant Measures For One Sided Dyck
Shift}\label{dyck1}

In this section we consider the one sided Dyck shift. We prove the following result:\\
\begin{theorem}\label{dyck1_ue}
The tail relation of the one sided Dyck shift is uniquely ergodic.
Furthermore, there exists a unique topologically $\sigma$-finite
tail-invariant measure on the one sided Dyck shift (up to
multiplication by a positive real number). \footnote{A measure
$\mu$ on topological space $X$ \emph{is topologically
$\sigma$-finite} if there is a countable cover of $X \setminus N$
by open sets with finite $\mu$-measure, where $N$ is a $\mu$-null
set.}
\end{theorem}

This theorem is a direct conclusion of lemmas \ref{G_0}, \ref{G_+}
and lemma \ref{G_-} below.
\begin{lemma}\label{top_trans_dyck}
The tail relation of the one sided $m$-Dyck is topologically
transitive.
\end{lemma}
\begin{proof}
Let $y = (y_n) \ \in \{\alpha_1,\ldots,\alpha_m\}^\mathbb{N}
\subset Y$. To prove the lemma, we will show that $\mathcal{T}(y)$
is dense in Y. Let $\omega = (\omega_1,\ldots,\omega_r) \in L(Y)$,
then $wy_{r+1}^{\infty} \in Y$. Thus, $[w] \cap \mathcal{T}(y)
\neq \emptyset$. This proves $\overline{\mathcal{T}(Y)}=Y$.
\end{proof}

If a tail-invariant measure $\mu$ on Y is topologically
$\sigma$-finite, $ \exists w \in L(Y) \mbox{ s.t. }
0<\mu([w])<\infty$. A corollary of our main result is that
any such $\mu$ is a finite measure.\\

 Define the following tail-invariant decomposition
of the one-sided Dyck shift: \[ G_{+}=\{ y \in Y :
\lim_{i\to\infty}H_{i}(y)=+\infty\}
\]
\[ G_{-}=\{y \in Y :
\liminf_{i\to\infty}H_{i}(y)=-\infty\}\]
\[ G_{0}= \{ y \in Y : \liminf_{i\to\infty}H_{i}(y) \in
(-\infty,+\infty)\}\] Obviously, \( Y=G_{+}\uplus G_{-} \uplus
G_{0} \)

Let
\[W_n=W^m_{n}:=\{l \in \mathcal{L}(Y,n) : \; l \equiv 1 (mod M) \}\]
where $M$ is the syntactic monoid of the $m$-Dyck shift. $W_{n}$
is the set of balanced words of length n. Denote:
\[w_{n}=w^m_n:=|W^m_{n}|\]
Let
\[\widetilde{W}_n^m=\widetilde{W}_n = \{l \in \mathcal{L}(Y,n) : \; l=\alpha_i \tilde{l} \beta_i \; ,
  1 \le i \le m\; , \tilde{l} \in W_{n-2}  \}\]
$\widetilde{W}_n$ is  the set of Dyck words of length $n$. Denote
$\tilde{w}_{n}=\tilde{w}^m_n:=|\widetilde{W}^m_{n}|$. Obviously,
$\tilde{w}_n \leq w_n$
\begin{lemma}\label{dyck_num_bal}
 \[w^m_{2k}=\frac{\left( \begin{array}{cc}2k \\ k
\end{array} \right)}{k+1}m^k\]
\end{lemma}
\begin{proof}
First, we note that \(w^m_{2k}=m^k w^1_{2k}\). This follows from
the fact that given $a \in W^1_{2k}$ one can independently choose
the "type" of each pair of brackets in order to create distinct
elements in $W^m_{2k}$, and every element of $W^m_{2k}$ can be
created this way. This describes a $m^k$ to one mapping $W^m_{2k}
\mapsto W^1_{2k}$.\\
All that remains is to prove
\[w^1_{2k}=\frac{\left( \begin{array}{cc}2k \\ k \end{array}
\right)}{k+1}\] This is sometimes called \emph{the ballot
problem}. An elementary proof of this can be found in pages 69-73
of \cite{FL}.
\end{proof}

\begin{lemma}\label{G_0}
 There are no topologically $\sigma$-finite
tail invariant measures,  giving $G_0$ positive measure.
\end{lemma}
\begin{proof}
Suppose $\mu$ is a tail invariant measure s.t. \( 0 <
\mu(G_{0}\cap [v])< \infty \) and \( |v|=l \). Without loss of
generality, we can assume \( \mu(G_{0} \cap [v] )=1 \) . Let
$R_{n}$ be the subset of $X$ consisting of points which are
balanced from time $n$ onwards:
\[R_n = \{y \in Y : \forall l \geq n \  H_l(y) \geq H_n(y),\  \liminf_{i\rightarrow \infty
}H_i(y)=H_n(y)\}\]

 We write the following decomposition of $[v] \cap R_n$, according to the
 first Dyck word following $v$:
\[\forall n \geq l , \; \
[v] \cap R_{n}=\biguplus_{k}\biguplus_{w\in
\widetilde{W}_{k}}(T^{-n}[w]\cap R_{n}) \cap [v] =
\biguplus_{k}\biguplus_{w\in \widetilde{W}_{k}}(T^{-n}[w]\cap
R_{n+k}) \cap [v] \] We further decompose each of these sets:
\[[v] \cap T^{-n}[w]\cap R_{n+k}=\biguplus_{a \in \mathcal{L}_n}(T^{-n}[w]\cap R_{n+k}
\cap [a]) \cap [v]= \biguplus_{a \in \mathcal{L}_n ,\\
a_{1}^l=v}[aw]\cap (R_{n+k}) \] We note that
\[R_{n+k} \cap [v] = \biguplus_{b \in
\mathcal{L}(Y,n+k) ,\\ a^l_1 = v} R_{n+k} \cap [b] \]
 and for every $a,b \in \mathcal{L}_{n+k}$, $\mu(R_{n+k} \cap [a])=\mu(R_{n+k}
\cap [b])$ because $\mu$ is $\mathcal{T}(Y)$-invariant.\\
 Let $r_{n}=\mu(R_{n} \cap [v])$, $a_{n}=|\{ \alpha \in \mathcal{L}_{n}:
\alpha_1^l=v\}|$,$r_{\infty}=\sup_{n}r_{n}$, then:
\[ r_{n}=\sum_{k}\frac{a_{n}\tilde{w}_{k}}{a_{n+k}}r_{n+k}\]
so:
\[r_{\infty}\leq
{\sup_{n \ge 0}(\sum_{k}\frac{a_{n}w_{k}}{a_{n+k}}})r_{\infty}\]
Since $r_n\le \mu([v]\cap G_0)$ and  $G_0=\bigcup_n R_n$, we have
$0 < r_{\infty} <\infty$ .We obtain:
\[1 \leq \sup_{n \ge 0}{\sum_{k}\frac{a_{n}w_{k}}{a_{n+k}}}\]
 Since for any \(u \in \mathcal{L}_n\)
and  any \( 1 \leq j \leq m \), \(u \alpha_j \in
\mathcal{L}_{n+1}\) and \(\exists 1 \leq j \leq m \; u\beta_{j}
\in \mathcal{L}_{n+1} \), we get the inequality  \(
\frac{a_{n+1}}{a_{n}} \geq m+1\). This proves \(
\frac{a_{n}}{a_{n+2k}} \leq \frac{1}{(m+1)^{2k}}\). Also,
 \[w_{2k}=\frac{\left( \begin{array}{cc}2k \\ k
\end{array} \right)}{k+1}m^k\]
from this follows:
\[\sum_{k=1}^{\infty}\frac{a_{n}w_{k}}{a_{n+k}} \leq
\sum_{k=1}^{\infty}\frac{ \left( \begin{array}{cc}2k \\ k
\end{array} \right)}{k+1}\left( \frac{m}{(m+1)^2} \right)^k \]

but:
\[\sum_{k=1}^{\infty}\frac{ \left( \begin{array}{cc}2k \\ k \end{array} \right)}{k+1}x^k =
x^{-1}\int_{0}^{x}\sum_{k=1}^{\infty}\left( \begin{array}{cc}2k \\
k \end{array} \right)t^{k} dt = \frac {1-\sqrt{1-4x}}{2x} -1
\]

so for $m>1$:
\[\sum_{k=1}^{\infty}\frac{a_{n}w_{k}}{a_{n+k}} \leq
\frac{(m+1)^2}{2m}\left( 1- \frac{m-1}{m+1} \right) -1  =
\frac{1}{m} \]

which implies:
\[\sup_{n \ge 0}{\sum_{k}\frac{a_{n}w_{k}}{a_{n+k}}} \le \frac{1}{m} < 1\]

This gives us a contradiction to our assumption of the existence
of such a measure $\mu$.
\end{proof}

\begin{lemma}\label{G_+}
 There are no topologically $\sigma$-finite
tail invariant measures, giving $G_+$ positive measure.
\end{lemma}
\begin{proof}
Assume there exist a tail invariant measure $\mu$  such that
$0<\mu(G_+ \cap [v])<\infty$. Since $G_{+}$ is a tail invariant
subset, we can assume $\mu(G_+^{c})=0$ by taking $\mu'(A)=\mu(A
\cap G_{+})$. Let:
\[\tilde{F}_{n}= \bigcap_{k > n}\{H_{k}(y)>H_{n}(y)\}\]
and:
\[F_{n}= \tilde{F}_n \setminus \bigcup_{j=1}^{n-1}\tilde{F}_j\]
$F_{n}$ is the sets of points which have the first $\alpha$ which
is unmatched at coordinate $n$. By definition, $$G_{+} \subset
\bigcup_{n
>0} F_{n}$$ so for some $n$ we must have $0< \mu(F_{n} \cap [v])<\infty$.
$$ F_n = \bigcup_w\bigcup_{i=1}^{m}(F_n \cap [w\alpha_i])$$
where the union is over all $w\in \mathcal{L}_{n-1}$ with
$\hat{\alpha}(w)=0$ . For any $K \in \mathbb{N}$ we have that:
$$(F_n \cap [w\alpha_i]) = \biguplus_b (F_n \cap [w\alpha_i b])$$
where this time the union is over all $b \in \mathcal{L}_{K}$ such
that $\hat{\beta}(b) = 0$. The reason there should be no unmatched
$\beta$'s in $b$ is so they will not match the $\alpha_i$ at
coordinate $n$. We denote the set of such $b$'s by $U_{K}$.
Suppose such $b$ has $\hat{\alpha}(b)=j$ with $j>0$. Denote by
$\xi (b,t)$, $0<t<j-1$ (which also depends on $w$), the word
obtained from $b$ by replacing the leftmost unmatched $\alpha_s$
by $\beta_i$ (so as to match the unmatched $\alpha_i$) and
replacing the next $t$ leftmost unmatched $\alpha_s$ with
$\beta_{s}$. It follows from the construction that for any $y \in
Y$, if $w\alpha_i b y \in Y$ then
$w\alpha_i\xi(b,t) y \in Y$.\\
 This shows there is a tail holonomy
$\pi:[w \alpha_i b] \to \pi[w\alpha_i b]\subset [w \alpha_i
\xi(b,t)]$, so $\mu([w\alpha_i b]) \le \mu([w\alpha_i\xi(b,t)])$.
For $b_1, b_2 \in U_{K}$, if $\xi(b_1,t) = \xi(b_2,t)$, this
implies that $b_1$ and $b_2$ can differ only where the first
unmatched $\alpha$ is located - so the maps $\xi(.,t)$ are $m$ to
$1$. Let
$$C(K,n,j):=\bigcap_{N>K} \{H_{N}(y) >j+H_{n+1}(y)\}$$
By definition, $G_{+} \subset \bigcup_{K}C(K,n,j)$. From our
assumption that $\mu$ is supported on $G_{+}$, it follows that
$\mu((\bigcup_{K}C(K,n,j))^{c}) =0$. Since $C(K,n,j)$ is an
increasing sequence of sets, there exist $K_0$ such that

$$\mu(C(K_0,n,j) \cap [w\alpha_i] \cap F_n ) > (1- \frac{1}{j})\mu([w\alpha_i] \cap F_n)$$

define:
$$U(K_0,j)=\{b \in \mathcal{L}_{K_0}:\; \hat{\alpha}(b) > j,\;\hat{\beta}(b)=0\}$$
 Note that if $t_1 \ne t_2$, then $\xi(b_1,t_1)
\ne \xi(b_2,t_2)$, because they have different number of unmatched
$\beta$'s. We have:
$$ C(K_0,n,j) \cap [w\alpha_i] \cap F_n = \biguplus_{b\in U(K_0,j)} ( [w\alpha_i b]) \cap C(K_0,n,j) \cap F_n$$

For the above $K_0$, the following inequalities hold:
$$ \mu(F_n \cap
[w\alpha_i]) \le \frac{j}{j-1}\mu(F_n \cap [w\alpha_i]\cap
C(K_0,n,j)) = $$
$$ \sum_{b \in U(K_0,j)}\mu([w\alpha_i b]) \cap C(K_0,n,j) \cap F_n) \le$$
$$ \sum_{b \in U(K_0,j)}\mu([w\alpha_i b]) \cap C(K_0,n,j)) \le$$
Because $\xi(.,t)$ are $m$ to 1:
$$ \le m \sum_{b \in U(K_0,j)}\mu([w\alpha_i \xi(b,t)] ) $$

We average this in equality over $t$:
$$ \mu(F_n \cap [w\alpha_i]) \le \frac{m}{j-1} \sum_{t=1}^{j-1}\sum_{b \in U(K_0,j)}\mu([w\alpha_i \xi(b,t)]
)$$
 because:
$$ [w\alpha_i]  \supseteq \biguplus_{t=1}^{j-1} \biguplus_{b\in U(K_0,j)} ( [w\alpha_i \xi(b,t)])$$
we obtain:
$$ \mu(F_n \cap [w\alpha_i]) \le \frac{m}{j-1} \mu([w\alpha_i])$$

We assume $ \mu([w\alpha_i]) \le \mu([v])=\mu(G_+ \cap [v])
<\infty$, so taking $j \rightarrow \infty$ we obtain that $\mu(F_n
\cap [w\alpha_i]) =0$, and since $G_+ \cap [v]$ is a countable
union of such sets we conclude that $\mu(G_+ \cap [v])=0$.
\end{proof}
 We conclude that every tail invariant measure of the Dyck shift is supported by
\[ G_{-}=\{ y \in X : \liminf_{i \rightarrow \infty}H_{i}(y)=-\infty\} \]
To prove unique ergodicity, we need the following:
\begin{lemma}\label{G_-}
There exists a unique tail-invariant probability measure \( \mu \)
on Y such that \( \mu(G_{-})=1 \). Furthermore, for any
topologically $\sigma$-finite tail-invariant measure $\mu'$ on
$G_{-}$, $\mu'=c \mu$ for some positive number $c$.
\end{lemma}
\begin{proof}
Let \( \Theta = \{\beta_1,...,\beta_m,\alpha\}^\mathbb{N}\).
Define \( \widetilde{H}_{0}(x)=0 \), \( \widetilde{H}_i(x) =
\sum_{j=1}^{i}(-\sum_{k=1}^{m}\delta_{x_j,\beta_k}+\delta_{x_j,\alpha})
\), $x\in \Theta$. Denote \[ \Theta_{-} = \{x \in \Theta :
\liminf_{i \to \infty} \widetilde{H}_{i}(x) = - \infty \} \]
 We will use a one-to-one Borel
mapping of \(G_{-} \) on to \( \Theta_{-} \), introduced in
\cite{WK74}. The map is defined is follows:
\[ g_-: G_{-}
\rightarrow \Theta_{-}
\]
\[ g_-(y)_{i}= \left\{
\begin{array}{ll}
\alpha & y_{i} \in \{ \alpha_{1},\ldots,\alpha_{m} \} \\
\beta_j & y_{i}= \beta_j
\end{array}
\right.
\]
$g_-$ is a bijection, and for any \( y_1,y_2 \in G_{-} \) \(
(y_1,y_2) \in \mathcal{T}(Y) \Leftrightarrow (g_-(y_1),g_-(y_2))
\in \mathcal{T}(\Theta) \). Let $p$ be the symmetric Bernoulli
measure on \( \Omega  \) satisfying \(
p([\omega_1,\ldots,\omega_n])=(\frac{1}{m+1})^n \). by the law of
large numbers \( p(\Theta_{-})=1 \), and therefor \( p \circ
g_{-}(G_{-})=1 \). So \( p \circ g_{-} \) is a tail invariant
probability measure on $Y$ supported by $G_{-}$. Suppose \( \mu \)
is a tail invariant probability measure on $Y$ s.t. \( \mu
(G_{-})=1 \). \( \mu \) can be transported by $g_-$ to a tail
invariant probability measure $q$ on \( \Theta \) (supported by \(
\Theta_{-} \)). Since \( \Theta \) is a full-shift, the uniqueness
of $\mathcal{T}(\Theta)$-invariant topologically $\sigma$-finite
measure follows immediately from the fact that all cylinders of
the same length have equal measure. This proves the uniqueness of
a tail-invariant topologically $\sigma$-finite measure on $G_{-}$.
\end{proof}

\section{Two Sided Dyck Shift}
\label{dyck2}
\subsection{Maximal Entropy Implies Double-Tail
Invariance}\label{mu_+_mu_-} In \cite{WK74} it was demonstrated
that the Dyck shift has two ergodic shift invariant probabilities
with entropy equal to the topological entropy. Such probabilities
are called \emph{equilibrium states}. In this section we show that
both of
these probabilities are also double-tail invariant.\\
 We introduce the following sets, which are mutually disjoint and
 are double-tail invariant. For $s,t \in \{\{+\infty\},\{-\infty\}, \mathbb
R\}$ we define:
\[B^s_t=\{x \in X : \liminf_{i \rightarrow + \infty}H_i(x) \in s ,
    \liminf_{i \rightarrow - \infty}H_i(x) \in t \}\]
let
\[\Omega^+_-=\{x \in \{\alpha_1 \ldots \alpha_m,\beta\}^{\mathbb{Z}} : \liminf_{i \rightarrow + \infty}\widehat{H}_i(x)=+\infty ,
 \liminf_{i \rightarrow - \infty}\widehat{H}_i(x)=-\infty \}\]
 and
 \[\Theta^-_+=\{x \in \{\beta_1 \ldots \beta_m,\alpha\}^{\mathbb{Z}} : \liminf_{i \rightarrow + \infty}\widetilde{H}_i(x)=-\infty ,
 \liminf_{i \rightarrow - \infty}\widetilde{H}_i(x)=+\infty \}\]
 Where $\widehat{H}$ and $\widetilde{H}$ are defined on
 $\{\alpha_1 \ldots \alpha_m,\beta\}^{\mathbb{Z}}$  and $\{\beta_1 \ldots \beta_m,\alpha\}^{\mathbb{Z}}$ respectively,
 as in formula (\ref{H_def}).

Define:
 \[g_+:B^{+\infty}_{-\infty} \mapsto \Omega^{+\infty}_{-\infty}\]
 \[ (g_+(y))_{i}= \left\{
\begin{array}{ll}
\alpha_j & y_{i}=\alpha_j  \\
\beta & y_{i} \in \{ \beta_{1},\ldots,\beta_{m}\}
\end{array}
\right.
\]
\[g_-:B^{-\infty}_{+\infty} \mapsto \Theta^{-\infty}_{+\infty}\]
 \[ (g_-(y))_{i}= \left\{
\begin{array}{ll}
\beta_j & y_{i}=\beta_j  \\
\alpha & y_{i} \in \{ \alpha_{1},\ldots,\alpha_{m}\}
\end{array}
\right.
\]
$g_+$ is a Borel bijection from $B^{+\infty}_{-\infty}$ to
$\Omega^{+\infty}_{-\infty}$ and $g_-$ is a Borel bijection of the
appropriate sets. The definitions of $g_+$ and $g_-$ can also be
extended to functions $g_+:B^{\mathbb{R}}_{-\infty} \mapsto
\Omega^{\mathbb{R}}_{-\infty}$ and $g_-:B^{-\infty}_{\mathbb{R}}
\mapsto \Theta^{-\infty}_{\mathbb{R}}$, which are also Borel
bijections.

\begin{lemma}\label{g_+_g_-_isomorphism}
$g_+:B^{+\infty}_{-\infty} \mapsto \Omega^{+\infty}_{-\infty}$,
$g_-:B^{-\infty}_{+\infty} \mapsto \Theta^{-\infty}_{+\infty}$,
$g_+:B^{\mathbb{R}}_{-\infty} \mapsto
\Omega^{\mathbb{R}}_{-\infty}$, $g_-:B^{-\infty}_{\mathbb{R}}
\mapsto \Theta^{-\infty}_{\mathbb{R}}$ are isomorphisms of the two
sided tail relations:
\[ (g_+ \times g_+) (\mathcal{T}_2(B^{+\infty}_{-\infty}))=
\mathcal{T}_2(\Omega^{+\infty}_{-\infty})\]
\[ (g_- \times g_-) (\mathcal{T}_2(B^{-\infty}_{+\infty}))=
\mathcal{T}_2(\Theta^{-\infty}_{+\infty})\]
\[ (g_+ \times g_+) (\mathcal{T}_2(B^{\mathbb{R}}_{-\infty}))=
\mathcal{T}_2(\Omega^{\mathbb{R}}_{-\infty})\]
\[ (g_- \times g_-) (\mathcal{T}_2(B^{-\infty}_{\mathbb{R}}))=
\mathcal{T}_2(\Theta^{-\infty}_{\mathbb{R}})\]

\end{lemma}

\begin{proof}
We prove the result for $g_+:B^{+\infty}_{-\infty} \mapsto
\Omega^{+\infty}_{-\infty}$, the other results are proved in the same manner.\\
$ (g_+ \times g_+) (\mathcal{T}_2(B^{+\infty}_{-\infty})) \subset
\mathcal{T}_2(\Omega^{+\infty}_{-\infty})$ is trivial, so we show
the other inclusion. Suppose $(g_+(x),g_+(y)) \in
\mathcal{T}_2(\Omega^{+\infty}_{-\infty})$. Let $n_0 \geq 0$ be
such that $g_+(x)_{[-n_0,n_0]^c}=g_+(y)_{[-n_0,n_0]^c}$.\\
Let \[r(i,x)=max\{j<i: H_j(x)=H_i(x)\}\] Clearly,
$r(i_1,x)=r(i_2,x)$ is impossible
 for $i_1 \neq i_2$.
 Since
\[\liminf_{n\rightarrow +\infty}H_n(x),\liminf_{n \rightarrow
+\infty}H_n(y)>-\infty\] there exists $c$ such that for some large
$N$, $H_i(x)>c$ for every $i>N$. Since $\liminf_{n \rightarrow
-\infty}H_n(x) = \liminf_{n \rightarrow -\infty}H_n(y)= -\infty$ ,
it follows that there exist some $i_0<N$ such that $H_{i_0}(x)=c$,
so for every $i>N$, $r(i,x)>i_0$. The same argument applies for
$y$. Since $(r(i,x))_{i>N}$ and $((r(i,y))_{i>N}$ are both
injective sequences of integers, bounded from below, it follows
that \[\lim_{n \rightarrow +\infty}r(n,x)=\lim_{n \rightarrow
+\infty}r(n,y)=+\infty\]
 Note that for $n_1,n_2> n_0$,
\[\widehat{H}_{n_1}(g_{+}(x))-\widehat{H}_{n_2}((g_{+}(x))=\widehat{H}_{n_1}(g_{+}(y))-\widehat{H}_{n_2}((g_{+}(y))\]
 so for all large $n$ enough so that $r(n,x)>n_0$,$r(n,y)>n_0$, there are exactly two cases:
 \begin{enumerate}
    \item $g_+(x)_n=g_+(y)_n=\beta$,in which case $r(n,x)=r(n,y)$ and $x_{r(n,x)}=y_{r(n,y)}$, so $x_n=y_n$
    \item $g_+(x)=g_+(y)=\alpha_i$ for $1< i < m$, and then $x_n= y_n=\alpha_i$
 \end{enumerate}
 Obviously, for $n< -n_0$, $x_n=y_n$. This proves $(x,y) \in
 \mathcal{T}_2(B^{+\infty}_{-\infty})$.
\end{proof}

\begin{lemma}\label{mu_+_mu_-_unique}
There exists a unique $\mathcal{T}_2$-invariant probability of $X$
supported by $B^{+\infty}_{-\infty}$, and a unique
$\mathcal{T}_2$-invariant probability of $X$ supported by
$B^{-\infty}_{+\infty}$. There are no $\mathcal{T}_2$-invariant
probabilities on $B^{\mathbb{R}}_{-\infty}$ and
$B^{-\infty}_{\mathbb{R}}$.
\end{lemma}

\begin{proof}
The symmetric product measure $p$ on $\Omega$ assigns probability
one to $\Omega^+_-$ Transporting the product measure on $\Omega$
by means of $g_+^{-1}$ to $B^{+\infty}_{-\infty}$ yields a tail
invariant
probability measure on $X$, by the previous lemma.\\
On the other hand, any tail invariant probability on $X$ supported
by $B^{+\infty}_{-\infty} \cup B^{\mathbb{R}}_{-\infty}$ can be
transported to a tail invariant probability on $\Omega$ by $g_+$.
This is an injective correspondence, so by the uniqueness of
double-tail invariant probability on $\Omega$, we conclude the
uniqueness of double-tail invariant probability on
$B^{+\infty}_{-\infty} \cup B^{\mathbb{R}}_{-\infty}$. In
particular, this also proves that no double-tail invariant
probability on $B^{\mathbb{R}}_{-\infty}$ exist. We obtain the
results for $B^{-\infty}_{+\infty}$ and $B^{-\infty}_{\mathbb{R}}$
symmetrically.
\end{proof}

\subsection{A Third Double-tail Invariant
Probability}\label{third_prob}
 For $z \in \{0,1\}^\mathbb{Z}$, we define:
\[ \tilde{H}_i(z)= \left\{ \begin{array}{ll}
  \sum_{j=0}^{i-1}
(\delta_{1,z_j}-\delta_{0,z_j}) & \mbox{if
$i>0$} \\
\sum_{j=i}^{-1}(\delta_{0,z_j}-\delta_{1,z_j})
& \mbox{if $i<0$} \\
0 & \mbox{if $i=0$} \\
\end{array}
\right.
\]
 Let \[S_{-\infty}^{-\infty}= \{z \in \{0,1\}^\mathbb{Z} : \;
\inf_{n \geq 0}\tilde{H}_n(z)=-\infty \; ,  \inf_{n <0
}\tilde{H}_n(z)=-\infty \;\}\] Let us define a Borel function $F:
S_{-\infty}^{-\infty} \times \{1,\ldots,m\}^\mathbb{Z}
\mapsto \Sigma^\mathbb{Z}$:\\

Let
\[F(z,a)_n= \left\{
\begin{array}{ll} \alpha_j & \mbox{if $z_n=1$, $a_{\gamma_n(z)}=j$} \\
\beta_j & \mbox{if $z_n=0$, $k=\varepsilon_n(z)$, and
$a_{\gamma_k(z)}=j$}
\end{array}
\right.
\]
where,
\[\gamma_k(z) = \left\{ \begin{array}{ll}
 \sum_{i=0}^{k}z_i & k \geq 0\\
 -\sum_{i=k}^{-1}z_i  & k<0\\
\end{array}
\right.\]
\[\varepsilon_n(z) = \max\{l < n : \;
\tilde{H}_l(z) \leq\tilde{H}_{n+1}(z)\}\]
 Since $\liminf_{n \rightarrow -\infty}\tilde{H}_n(z)=-\infty$
for $z \in S_{-\infty}^{-\infty}$, $F$ is well defined.\\
\begin{lemma}
For every $z \in S_{-\infty}^{-\infty}$, $a \in
\{1,\ldots,m\}^\mathbb{Z}$, $F(z,a) \in X$.
\end{lemma}
\begin{proof}
Suppose $x=F(z,a) \not\in X$, then there exist $n,n' \in
\mathbb{Z}$, $n<n'$, such that $x_n=\alpha_i$, $x_{n'}=\beta_j$
with $i \neq j$ and $n=\max\{l < n' : \; H_l(x)=H_{n'+1}(x)\}$.
But in that case, $n=\varepsilon_{n'}(z)$, so
$i=j=a_{\gamma_n(z)}$.
\end{proof}
 Let $\mu_1$ be the symmetric product measure on $\{0,1\}^\mathbb{Z}$, and
$\mu_2$ the symmetric product measure on
$\{1,\ldots,m\}^\mathbb{Z}$.
\begin{lemma}
$\mu_1(S_{-\infty}^{-\infty})=1$
\end{lemma}
\begin{proof}
This follows from recurrence and ergodicity of the simple random
walk on $\mathbb{Z}$.
\end{proof}
We define: $\tilde{\mu}=(\mu_1 \times \mu_2)\circ F^{-1}$. Since
$F^{-1}(B_{-\infty}^{-\infty})=S^{-\infty}_{-\infty} \times
\{1,\ldots,m\}^\mathbb{Z}$ it follows
that $\tilde{\mu}(B_{-\infty}^{-\infty})=1$.\\

Let us also define a Borel mapping $z:B_{-\infty}^{-\infty}
\mapsto S_{-\infty}^{-\infty}$:
\[z(x)_n= \left\{ \begin{array}{ll} 1 & x_n \in \{\alpha_1,\ldots,\alpha_m\} \\
0 & x_n \in \{\beta_1,\ldots,\beta_m\} \end{array} \right. \]

The following lemma gives an explicit formula for the
$\tilde{\mu}$ probability of a cylinder:
\begin{lemma}\label{tilde_mu_explicit}
Let $w \in L(X)$. If the number of matched $\alpha$'s in $w$ is
$n_1$ and the number of unmatched $\alpha$'s and $\beta$'s is
$n_2$ ($2n_1+n_2=|w|$) then
$\tilde{\mu}([w]_k)=m^{-(n_1+n_2)}(\frac{1}{2})^{|w|}$.
\end{lemma}
\begin{proof}
Denote by $f_1,\ldots,f_{n_1}$ the locations of matched $\alpha$'s
in $w$. Denote by $g_1,\ldots,g_{n'_2}$ the locations of unmatched
$\alpha$'s in $w$. Denote by $h_1,\ldots,h_{n''_2}$ the locations
of unmatched $\beta$'s in $w$. We have $n'_2+n''_2=n_2$. For
$\vec{r} \in \mathbb{Z}^{n_1}$, $\vec{s} \in \mathbb{Z}^{n'_2}$,
$\vec{t} \in \mathbb{Z}^{n''_2}$, define:
\[A_{\vec{r}}=\{z :\; \gamma_{k+f_l}(z)=r_l \; 1\leq l\leq n_1 \}\]
\[B_{\vec{s}}=\{z :\; \gamma_{k+g_l}(z)=s_l \; 1\leq l\leq n'_2 \}\]
\[C_{\vec{t}}=\{z :\; \gamma_{\varepsilon_l}(z)=t_l \; \varepsilon_l=\varepsilon_{k+h_l}(z) 1\leq l\leq n''_2 \}\]
Informally, $A_{\vec{r}},B_{\vec{s}},C_{\vec{t}}$ determine the
locations in the sequence $a \in \{1,\ldots,m\}^{\mathbb{Z}}$
involved in selecting the types of $\alpha$'s and $\beta$'s within
the coordinates $k,\ldots,k+|w|$. Now we define:
\[Z = \{z \in S_{-\infty}^{-\infty} :\; z_{i+k} = z(w)_i \mbox{ for } 0 \leq i
\leq |w| \}\]
\[A'_{\vec{r}}=\{a \in \{1,\ldots,m\}^{\mathbb{Z}}:\; a_{r_l}=j \mbox{ if } w_{f_l}=\alpha_j\}\]
\[B'_{\vec{s}}=\{a \in \{1,\ldots,m\}^{\mathbb{Z}}:\; a_{s_l}=j \mbox{ if } w_{g_l}=\alpha_j\}\]
\[C'_{\vec{t}}=\{a \in \{1,\ldots,m\}^{\mathbb{Z}}:\; a_{t_l}=j \mbox{ if } w_{h_l}=\beta_j\}\]
With the above definitions we can write:
\begin{equation}
F^{-1}([w]_k)= Z \times \{1,\ldots,m\}^{\mathbb{Z}}\cap
\bigcup_{\vec{s},\vec{t},\vec{r}}((A_{\vec{r}}\times
A'_{\vec{r}})\cap(B_{\vec{s}}\times B'_{\vec{s}}) \cap
(C_{\vec{t}}\times C'_{\vec{t}}))\
\end{equation}
Where the union of $\vec{r},\vec{s},\vec{t}$ ranges over all
vectors such that the set of numbers appearing in their
coordinates are pairwise disjoint. This is a union of disjoint
sets. Thus:

\[\tilde{\mu}([w]_k)=\sum_{\vec{s},\vec{t},\vec{r}}(\mu_1 \times
\mu_2) ((Z\cap A_{\vec{r}} \cap B_{\vec{s}} \cap C_{\vec{t}})
\times (A'_{\vec{r}}\cap B'_{\vec{s}} \cap C'_{\vec{t}}))\]
\begin{equation} \label{mu_1_prod_mu_2}
\tilde{\mu}([w]_k)=\sum_{\vec{s},\vec{t},\vec{r}}\mu_1(Z\cap
A_{\vec{r}} \cap B_{\vec{s}} \cap
C_{\vec{t}})\mu_2(A'_{\vec{r}}\cap B'_{\vec{s}} \cap C'_{\vec{t}})
\end{equation}
 Now notice that for every $\vec{r},\vec{s},\vec{t}$ in the sum,
\[\mu_2(A'_{\vec{r}}\cap B'_{\vec{s}} \cap
C'_{\vec{t}})=m^{-(n_1+n'_2+n''_2)}=m^{-(n_1+n_2)}\] Also note
that $Z = \biguplus_{\vec{s},\vec{t},\vec{r}}(Z\cap A_{\vec{r}}
\cap B_{\vec{s}} \cap C_{\vec{t}})$, so
$\mu_1(Z)=\sum_{\vec{s},\vec{t},\vec{r}}\mu_1(Z\cap A_{\vec{r}}
\cap B_{\vec{s}} \cap C_{\vec{t}})$. Thus, equation
\ref{mu_1_prod_mu_2} can be simplified as follows:
\[\tilde{\mu}([w]_k)=\sum_{\vec{s},\vec{t},\vec{r}}\mu_1(Z\cap
A_{\vec{r}} \cap B_{\vec{s}} \cap
C_{\vec{t}})m^{-(n_1+n_2)}=\mu_1(Z)m^{-(n_1+n_2)}=(\frac{1}{2})^{|w|}m^{-(n_1+n_2)}\]
\end{proof}

\begin{theorem}\label{tilde_mu}
$\tilde{\mu}$ is a $\mathcal{T}_2$-invariant probability.
\end{theorem}
Our method of proving this is as follows:
 We define a countable set of $\mathcal{T}_2$-holonomies
 \[\mathcal{H}=\{g_{w,w',n}: \; n \in \mathbb{Z}, \; w,w' \in L(X) \; |w|=|w'|,\; w \equiv w' \mbox{(mod $M$)}, \; \}\]
 By proposition \ref{g_w_w'} below, we see that $\tilde{\mu}$ is
 invariant under $\mathcal{H}$. Then we prove that $\mathcal{H}$
 generates $\mathcal{T}_2$, up to a $\tilde{\mu}$-null set (proposition
 \ref{H_generate_mod_mu} bellow). This will complete the proof.

\begin{lemma}\label{g_is_good}
Suppose $w,w' \in \mathcal{L}(X,n)$ with $w \equiv w'$ (mod $M$).
If $x,y \in \Sigma^{\mathbb{Z}}$ such that $x_{[k-n,k]}=w$,
 $y_{[k-n,k]}=w'$ and $x_{[k-n,k]^c}=y_{[k-n,k]^c}$, then
\[x \in X \Leftrightarrow y \in X\]
\end{lemma}
\begin{proof}
Suppose $x \in X$. We have to show that for every $j
> n$ $y_{[k-j,j]} \not\equiv 0$ (mod $M$). Writing
$x_{[k-j,j]}=swt$ , we have $y_{[k-j,j]}=sw't$ and since $ w\equiv
w'$ (mod $M$), $sw't \equiv swt \not\equiv 0$ (mod $M$). This
shows $y \in X$. By replacing the roles of $y$ and $x$ we get: $y
\in X \Rightarrow x \in X$.
\end{proof}

 Let $w,w' \in \mathcal{L}(X,n)$ with $w \equiv w'$ (mod $M$) and $k
\in \mathbb{Z}$. By lemma \ref{g_is_good} we can define
$g_{w,w',k}:[w]_k \mapsto [w']_k$ to be the Borel function that
changes the $n$ coordinates starting at $k$ from $w$ to $w'$.
\[g_{w,w',k}(\ldots,x_{k-1},w_0,\dots,w_{n-1},x_{k+n},\ldots)=(\ldots,x_{k-1},w'_0,\ldots,w'_{n-1},x_{k+n}\ldots)\]

\begin{proposition}\label{w_w'}
If $w \equiv w' (\mbox{mod } M)$ , $|w|=|w'|$, and $k \in
\mathbb{Z}$, then $\tilde{\mu}([w]_k)=\tilde{\mu}([w']_k)$.
\end{proposition}
\begin{proof}
By lemma \ref{tilde_mu_explicit},
$\tilde{\mu}([w]_k)=m^{-(n_1+n_2)}(\frac{1}{2})^{|w|}$. Since the
number of paired $\alpha$ in $w'$ is also $n_1$, we get that
$\tilde{\mu}([w]_k)=\tilde{\mu}([w']_k)$.
\end{proof}
\begin{proposition}\label{g_w_w'}
If $w \equiv w' (\mbox{mod } M)$, $|w|=|w'|$, and $k \in
\mathbb{Z}$, then $\tilde{\mu}$ is $g_{w,w',k}$ invariant.
\end{proposition}
\begin{proof}
First note that if $w \equiv w' (\mbox{mod } M)$ then for every
$s,t \in L(X)$ $swt \equiv sw't (\mbox{mod } M)$. This fact, along
with proposition \ref{w_w'} shows that
$\tilde{\mu}(A)=\tilde{\mu}(g_{w,w',k}(A))$ for every cylinder set
A. Since the cylinder sets generate the Borel sets, this shows
$\tilde{\mu}$ is $g_{w,w',k}$-invariant.
\end{proof}
For $x \in B_{-\infty}^{-\infty}$, and $j >0$ define:
\[a_j(x)=\min \{ k >0 : \; H_{k+1}(x)=-j\}\]
\[b_j(x)=\max \{ k <0 : \; H_k(x)=-j\}\]
Note that for any $x \in B_{-\infty}^{-\infty}$, $(a_j(x))_{j \in
\mathbb{N}}$ is strictly increasing, and $(b_j(x))_{j \in
\mathbb{N}}$ is strictly decreasing. Also note that $x_{a_j(x)}
\in \{\beta_1,\ldots,\beta_m\}$ and $x_{b_j(x)} \in
\{\alpha_1,\ldots,\alpha_m\}$, and if $x_{a_j(x)}=\beta_i$ then
$x_{b_j(x)}=\alpha_i$.
 Let $A_c^n =\{x \in B_{-\infty}^{-\infty}:\;x_{b_j(x)}=x_{b_{j+c}(x)}\; \forall j> n\}$.

\begin{lemma}\label{mu_A_c_0}
$\tilde{\mu}(A_c^n)=0$ for all $c \in \mathbb{Z}\setminus \{0\}$,
$n \geq 0$
\end{lemma}
\begin{proof}
For $z \in S_{-\infty}^{-\infty}$ define
\[\tilde{b}_j(z)=\max \{ k < 0: \; \tilde{H}_k(z)=j\}\]
For any  $x \in B_{-\infty}^{-\infty}$,
$\tilde{b}_j(z(x))=b_j(x)$. Now, for $J \subset \mathbb{N}$ with
$|J| < \infty$:
\[ \tilde{\mu}(\{x_{b_j(x)}=x_{b_{j+c}(x)} \; \mbox{for $j\in J$
}\})=\]
\[ (\mu_2\times\mu_1)(\{(a,z) : \; a_{l_{j,1}} = a_{l_{j,2}},\;
l_{j,1}=\tilde{b}_j(z)\; l_{j,2}=\tilde{b}_{j+c}(z) \; \mbox{for
$j \in J$} \})= (\frac{1}{m})^{|J|}\]
 This follows from the definition of
$\tilde{\mu}$ as the image of a product measure, and from the fact
that $(b_j(x))_{j \in \mathbb{N}}$ is strictly monotonic, so the
$l_{j,1}$'s are all distinct, and $l_{j,1} \neq l_{j,2}$ for $j
\in J$. Thus, $\tilde{\mu}(A_c^n)=0$.
\end{proof}
\begin{proposition}\label{H_generate_mod_mu}
 There exists a double-tail invariant set $X_0 \subset X$ with $\tilde{\mu}(X_0)=1$, such the countable set
of $\mathcal{T}_2$-holonomies
\[\mathcal{H}=\{g_{w,w',n}: \; n \in \mathbb{Z}, \; w,w' \in L(X) \; |w|=|w'|,\; w \equiv w' \mbox{(mod $M$)}, \; \}\]
generates $\mathcal{T}_2(X_0)$.
\end{proposition}

\begin{proof}
Let  $X_0 = B_{-\infty}^{-\infty}\setminus \bigcup_{n,m
> 0}\bigcup_{c \neq
0}T^{-m}A_c^n$. Since $\tilde{\mu}(B_{-\infty}^{-\infty})=1$, and
$\tilde{\mu}(A_c^n)=0$ for $c \neq 0$ by the previous lemma,
$\tilde{\mu}(X_0)=1$. Also, since $B_{-\infty}^{-\infty}$ and
$\bigcup_{n,m
> 0}\bigcup_{c \neq
0}T^{-m}A_c^n$ are $\mathcal{T}_2$-invariant sets, $X_0$ is
$\mathcal{T}_2$-invariant.
We show that $\mathcal{H}$ generates $\mathcal{T}_2(X_0)$.\\
Suppose $(x,y) \in \mathcal{T}_2(X_0)$. We must show that $y=
g(x)$ for some $g \in \mathcal{H}$. $\exists n \in \mathbb{N}$ so
that $x_{[-n,n]^c}=y_{-[n,n]^c}$.
 Let $w=x_{[-n,n]}$, $w'=y_{[-n,n]}$. Let $c=H(w)-H(w')$.\\
First assume $c \neq 0$. Let $x'=T^{-n}(x)$, $y'=T^{-n}(y)$. Then
$x'_{[0,2n]^c}=y'_{[0,2n]^c}$. For all $k > 2n$,
$H_k(x')=H_k(y')+c$. Therefore, $a_j(x')=a_{j+c}(y')$ for all
$j>2n+|c|$. Also, Since $x'_{[0,2n]^c}=y'_{[0,2n]^c}$,
$H_k(x')=H_k(y')$ for all $k<0$. So $b_j(x')=b_j(y')$ for all $j
>0$.\\
 For $j>2n+|c|$, denote $x'_{a_j(x')}=\beta_i$. Then
 $x'_{b_j(x')}=\alpha_i$. Also,
 $y'_{a_{j+c}(y')}=y'_{a_j(x')}=x'_{a_j(x')}=\beta_i$, so
 $y'_{b_{j+c}(y')}=\alpha_i$. Therefore,
 $x'_{b_{j+c}(x')}=y'_{b_{j+c}(y')}=\alpha_i$. We conclude
 that $x'_{b_j(x')}=x'_{b_{j+c}(x')}$ for all $j > 2n+|c|$.
  This proves that $x \in T^{-n}A_{c}^{2n+|c|}$, but we assumed $x
\in X_0$, so this is a
contradiction, so  $c=0$.\\
Therefore, for every $k_1<-n$ and $k_2>n$, we have:
\[H_{k_1}(x)-H_{k_2}(x)= H_{k_1}(y)-H_{k_2}(y)\]
 Let $N=\min \{k \geq n : \; H_{k+1}(x)<-2n\}$, and $N'=\max \{k<-n: \;
H_k(x)=H_{N(x)+1}(x)\}$. $N$ and $N'$ are well defined for $x \in
B_{-\infty}^{-\infty}$. We have that $H_{N+1}(x)-
H_{N'}(x)=H_{N+1}(y)- H_{N'}(y)=0$, and so $x_{[N',N]} \equiv
y_{[N',N]} \equiv 0$ (mod
$M$). Thus $y=g_{x_{[N',N]},y_{[N',N]},N'}(x)$.\\
\end{proof}

\begin{proposition}\label{tilde_mu_shift}
$\tilde{\mu}$ is a shift invariant probability.
\end{proposition}
\begin{proof}
Let $[w]_k$ be a cylinder set.By lemma \ref{tilde_mu_explicit}, we
have:
\[\tilde{\mu}([w]_k)=m^{-n_1+n_2}(\frac{1}{2})^{|w|}\]
and also:
\[\tilde{\mu}(T^{-1}[w]_k)=m^{-n_1+n_2}(\frac{1}{2})^{|w|}\]
So $\tilde{\mu}(A)=\tilde{\mu}(T^{-1}[A])$ for every Borel set
$A$.
\end{proof}
One could question whether proposition \ref{tilde_mu_shift}
follows immediately from the fact that the shift mapping is a
normalizer of the double-tail relation. We note that in general
double-tail invariant measures are not necessarily shift
invariant. To see this, consider a (finite) subshift consisting of
an orbit of a periodic point. For more elaborate examples of a
similar phenomenon  see \cite{BS96}, where it is shown that the
"generalized hard core model" has Gibbs measures which are not
shift-invariant.

\begin{proposition}\label{h_tilde_mu}
\[h_{\tilde{\mu}}(X,T)=\log(2)+\frac{1}{2}\log(m)\]
\end{proposition}
\begin{proof}
We have $h_{\tilde{\mu}}(X,T)=\lim_{n \rightarrow
\infty}h_{\tilde{\mu}}(x_0 | x_{-1},x_{-2},\ldots,x_{-n})$. Let
\[\varpi(a_1,\ldots,a_n)=\min\{H(a_1,\ldots,a_k)
: \; 0 \leq k \leq n\}\] By applying lemma
\ref{tilde_mu_explicit}, we get:
\[h_{\tilde{\mu}}(x_0 | x_{-1}=a_1,\ldots,x_{-n}=a_{n})=
\left\{ \begin{array}{ll} \log(2m) & \mbox{if
$\varpi(a_1,\ldots,a_n) \geq 0$} \\
\log(2)+\frac{1}{2}\log(m) & \mbox{if $\varpi(a_1,\ldots,a_n) <
0$}\end{array} \right.
\]
We have $h_{\tilde{\mu}}(x_0 |
x_{-1},x_{-2},\ldots,x_{-n})=\tilde{\mu}(\varpi(a_1,\ldots,a_n) <
0)(\log(2)+\frac{1}{2}\log(m))+\tilde{\mu}(\varpi(a_1,\ldots,a_n)
\geq 0)\log(2m)$. Since $\lim_{n \rightarrow
\infty}\tilde{\mu}(\varpi(a_1,\ldots,a_n) \geq 0)=0$, we have
$h_{\tilde{\mu}}(X,T)=\log(2)+\frac{1}{2}\log(m)$.
\end{proof}
For $m \geq 2$, $h_{\tilde{\mu}}(X,T)<h_{\mathit{top}}(X,T)$.
Thus, $\tilde{\mu}$ provides an example of a shift invariant
probability, which is also $\mathcal{T}_2$ invariant, yet has
entropy which is strictly less than the topological entropy, for
$m \geq 2$ .

\subsection{No other Double-Tail Invariant
Probabilities}\label{no_more}
 In this subsection we conclude that
apart from the two probabilities described in section
\ref{mu_+_mu_-} and the probability defined in section
\ref{third_prob}, there are no other ergodic double-tail
invariant probabilities for the Dyck shift.\\
By lemma \ref{mu_+_mu_-_unique} we know that there are no more
double-tail invariant probabilities  on the sets
$B^{+\infty}_{-\infty}$ and $B^{-\infty}_{+\infty}$. We also know
by the same lemma that there are no such probabilities on
$B^{\mathbb{R}}_{-\infty}$ and $B^{-\infty}_{\mathbb{R}}$.

Our next goal is to prove $\tilde{\mu}$ is unique on
$B^{-\infty}_{-\infty}$:\\
\begin{proposition}\label{Balance_prob_infty_infty}
Suppose $\nu$ is a $\mathcal{T}_2(B_{-\infty}^{-\infty})$
invariant probability. Then for every $w \equiv 1$ (mod $M$),
\[\nu([w]_t)= (\frac{1}{2\sqrt{m}})^{|w|}\]
\end{proposition}
\begin{proof}
Let $[w]_t$ be a balanced cylinder with $|w|=2n$. For $i<t$,
Denote:
\[M_{i,i+2N}=\{x \in X: \; x_i^{i+2N} \equiv 1 \mbox{( mod $M$)}\}\]
Since all balanced cylinders of the same length have equal $\nu$-
probability, we can calculate $\nu([w]_t \mid M_{i,i+2N})$ by
counting the number of balanced words of length $2N$, and the
number of such balanced words  with $w$ as a subword starting at
position $t-i$. By lemma \ref{dyck_num_bal}, the number of
balanced words of length $2N$ is
\[w^m_{2N}=\frac{\left( \begin{array}{cc}2N \\ N
\end{array} \right)}{N+1}m^N\]
The number balanced word of length $2N$ with $w$ as a subword
starting at position $t-i$ is $w^m_{2N-2n}$. Thus,
\[\nu([w]_t \mid M_{i,i+2N})=\frac{w^m_{2N-2n}}{w^m_{2N}}\]
It easily follows that:
\[\lim_{N \rightarrow \infty}\nu([w]_t \mid M_{i,i+2N})=
\lim_{N \rightarrow
\infty}\frac{w^m_{2N}}{w^m_{2N-2n}}=(\frac{1}{2\sqrt{m}})^{2n}\]
Since $\nu(B_{-\infty}^{-\infty})=1$, we have
\[\nu(\bigcap_{N_0 \in \mathbb{N}}\bigcup_{i \in -\mathbb{N}}
\bigcup_{N > N_0}M_{i,i+2N})=1\] For $N_0 > n$ define a random
variable $ \chi_{N_0}(x) := \min\{ N > N_0 : \; x \in \bigcup_{i
\in -\mathbb{N}}M_{i,i+2N} \}$. We have
\[\nu([w]_t)=\sum_{N>N_0}\nu(\chi_{N_0} = N)\nu([w]_t \mid
\chi_{N_0}=N) \rightarrow(\frac{1}{2\sqrt{m}})^{2n}\]
\end{proof}

\begin{proposition}
$\tilde{\mu}$ is the unique $\mathcal{T}_2$ invariant probability
on $B_{-\infty}^{-\infty}$.
\end{proposition}
\begin{proof}
Suppose $\nu$ is a $\mathcal{T}_2$ invariant probability on
$B_{-\infty}^{-\infty}$. By proposition
\ref{Balance_prob_infty_infty},
\begin{equation} \label{bal_cond}
 \forall w \equiv 1 \mbox{(mod $M$)} \; \nu([w])=
(\frac{1}{2\sqrt{m}})^{|w|}
\end{equation}
For $ a \in L(X)$, we say that $w \in L(X)$ is a \emph{minimal
balanced extension} of $a$, if the following conditions hold:
\begin{enumerate}
    \item There exist $l,r \in L(X)$ such that $w=lar$.
    \item $w \equiv 1$ (mod $M$)
    \item For every  $l'$ suffix of $l$ and $r'$ prefix of $r$,
    $l'ar' \equiv 1$ implies $l'ar' = w$.
\end{enumerate}
Since for every $a \in L(X)$,
 \[[a]_t =_{\nu} \biguplus\{[w]_s : \mbox{$w$ is a minimal balanced extension of $a$, with $(w_i)_{i=t-s}^{t-s+|w|}=a$}\}\]
 We have:
 \[\nu([a]_t)= \sum_{[w]_s}\nu([w]_s)=\sum_{[w]_s}\tilde{\mu}([w]_s)=\tilde{\mu}([a]_t)\]
 Where the sum ranges over minimal balanced extensions of $a$.
 This proves $\nu=\tilde{\mu}$.
 By theorem \ref{tilde_mu}, this proves $\tilde{\mu}$ is the unique double tail invariant
probability of $B_{-\infty}^{-\infty}$.
\end{proof}

Finally, we show that no other double-tail invariant
probabilities exist for the Dyck Shift.\\
Define:
 $\hat{p}:\Sigma^{\mathbb{Z}} \mapsto \Sigma^{\mathbb{N}}$ by $\hat{p}((x_n)_{n \in
\mathbb{Z}})=(x_n)_{n \in \mathbb{N}}$. This is a Borel mapping
that maps the two-sided Dyck shift $X$ onto the
one sided Dyck shift $Y \subset \Sigma^\mathbb{N}$.\\
 Let $K_0= \{ x \in X : H_i(x) \geq 0 ,\forall i<0 \}$, and $K_i=T^{-i}(K_0))$. Notice that
$B^s_t \subset \bigcup_{i=0}^{\infty}K_i$, for $s,t \in \{\{+\infty\},\mathbb{R}\}$. \\
\begin{lemma}\label{T_2_to_T}
If $A,B \subset Y$ are Borel sets, and $g:A \mapsto B$ is a
$\mathcal{T}(Y)$-holonomy, then there exists a
$\mathcal{T}_2(X)$-holonomy $\tilde{g}:(\hat{p}^{-1}(A)\cap K_0)
\mapsto (\hat{p}^{-1}(B)\cap K_0)$
\end{lemma}
\begin{proof}
We define $\tilde{g}:(\hat{p}^{-1}(A)\cap K_0) \mapsto
(\hat{p}^{-1}(B)\cap K_0)$ as follows:
\[\tilde{g}(x)_n = \left\{ \begin{array}{ll}
x_n & n<0 \\
g(\hat{p}(x))_n  & n\geq 0 \end{array} \right.\] We prove that
$\tilde{g}$ takes $\hat{p}^{-1}(A)\cap K_0$ into
$\hat{p}^{-1}(B)\cap K_0$. Let $x \in \hat{p}^{-1}(A)\cap K_0$.
 Since $x_n=\tilde{g}(x)_n$ for all $n <0$, we have $H_n(x)=
 H_n(\tilde{g}(x))$ for $n <0$. Because $x \in K_0$ we have
 $H_n(\tilde{g}(x)) \geq 0$ for $i < 0$. Let $y=\tilde{g}(x)$. Now we prove that  $y \in X$. Otherwise, there exist $n_1,n_2 \in \mathbb{Z}$, such that
$n_1=\min\{l<n_2:H_l(y)=H_{n_2+1}(y)\}$, and $y_{n_1}=\alpha_i$
$y_{n_2}=\beta_j$ with $i \neq j$. If $n_1,n2  <0$ then
$y_{n_1}=x_{n_1}$, $y_{n_2}=x_{n_2}$, so this contradicts the fact
that $x \in X$. If $n_1,n2  \geq 0$, then
$y_{n_1}=g(\hat{p}(x))_{n_1}$, $y_{n_2}=g(\hat{p}(x))_{n_2}$, so
this contradicts the fact that $g(\hat{p}(x)) \in Y$.\\
We remain with the case $n_1 < 0 \leq n_2$. We have $H_{n_1}(y)
\geq 0 = H_0(y)$, and $H_{n_2+1}(y)=H_{n_2}(y)-1$ (since
$y_{n_2}=\beta_j$).Also, $H_{n_2+1}(y) = H_{n_1}(y) \geq 0$. Since
$H_i(y)-H_{i+1}(y)= \pm 1$, there must be some $l >0$ such that
$H_l(y)=H_{n+1}(y)$. This contradicts the condition on $n_1,n2$.
By the definition of $\tilde{g}$,
$\hat{p}(\tilde{g}(x))=g(\hat{p}(x))$, so $\tilde{g}(x) \in
\hat{p}^{-1}(B)$. The fact that $g$ is one to one and onto
$(\hat{p}^{-1}(B)\cap K_0)$ follows from the fact that
\[\tilde{g}^{-1}(x)_n = \left\{ \begin{array}{ll}
x_n & n<0 \\
g^{-1}(\hat{p}(x))_n  & n\geq 0 \end{array} \right.\] To complete
the proof of  the lemma we must show that $(x,\tilde{g}(x)) \in
\mathcal{T}_2(X)$. Since $g$ is a $\mathcal{T}(Y)$-holonomy,
$\hat{p}(x)$ and $g(\hat{p}(x)$ only differ in a finite number of
(positive) coordinates. $x$ and $\tilde{g}(x)$ only differ in the
coordinates where $\hat{p}(x)$ and $g(\hat{p}(x))$ differ, which
is a finite set. So $(x,\tilde{g}(x)) \in \mathcal{T}_2(X)$
\end{proof}

\begin{lemma}
There are no $\mathcal{T}_2(X)$-invariant probability measures on
X supported by $B^s_t$, $s,t \in \{ \{+\infty\},\mathbb{R}\}$.
\end{lemma}
\begin{proof}
We first prove the result for $B^{\mathbb{R}}_t$,$t \in \{
\{+\infty\}, \mathbb{R}\}$. Recall that $K_i=\{x \in X:\; H_n(x)
\ge H_i(x), \forall n<i\}$. Notice that $B^{\mathbb{R}}_t \subset
\bigcup_{i=0}^{\infty}K_i$. \\
Suppose $\mu$ is a $\mathcal{T}_2(X)$-invariant probability
supported by $B^{\mathbb{R}}_t$,where $t \in \{ \{+\infty\},
\mathbb{R}\}$, then $\mu(K_i)>0$ for some $i \ge 0$.
Without loss of generality we can assume $\mu(K_0)>0$.\\
Define a probability $\breve{\mu}$ on $Y$ by the formula:
\[\breve{\mu}(A)=\frac{\mu(\hat{p}^{-1}(A)\cap K_0)}{\mu{K_0}}\]
By lemma \ref{T_2_to_T}, $\breve{\mu}$ is a $\mathcal{T}(Y)$
invariant probability. Also, since $\mu(B^{\mathbb{R}}_t)=1$,
\[\breve{\mu}(\{y \in Y: \liminf_{n \rightarrow + \infty}H_n(y) \in  \mathbb{R}\})=1\]
 Similarly, the existence of a $\mathcal{T}_2(X)$-invariant probability supported
by $B^{+\infty}_t$,where $t \in \{ \{+\infty\}, \mathbb{R}\}$
would result in a $\mathcal{T}(Y)$-invariant probability
$\breve{\mu}$ with
\[\breve{\mu}(\{y \in Y: \liminf_{n \rightarrow + \infty}H_n(y) =  +\infty \})=1\]
But in section \ref{dyck1} it was proved that the one sided Dyck
shift has a unique $\mathcal{T}$-invariant probability, supported
by
\[\{y \in Y: \liminf_{n \rightarrow + \infty}H_n(y) = -\infty\}\]

Which rules out the possibility that such $\breve{\mu}$ exists.

\end{proof}

\section*{Acknowledgements}
The Author thanks Jon Aaronson for his guidance during this work.
Also, the author thanks the anonymous referee for his helpful
comments.


\begin{thebibliography}{alpha}

\bibitem{ANS}
J. Aaronson, H. Nakada and O. Sarig, Exchangeable Measures For
Subshifts, http://www.arxiv.org/abs/math.DS/0406578


\bibitem{BH86}
F. Blanchard and G. Hansel, Systems Codes, Theortical Computer
Science 44 14-49, 1986.

\bibitem{BM}
R. Bowen and B. Marcus,  Unique ergodicity for horocycle
foliations, Israel J. Math. 26 no. 1, P. 43-67, 1977.

\bibitem{BS96}
R. Burton and  J. Steif, Some $2$-d symbolic dynamical systems:
entropy and mixing.  Ergodic theory of $Z\sp d$ actions , London
Math. Soc. Lecture Note Ser.,  Cambridge Univ. Press, no. 228, P.
297-305, 1996.


\bibitem{B05}
J. Buzzi, Subshifts of Quasi-Finite Type, Invert. Math. 159 P.
369-406, 2005.

\bibitem{FM}
J. Feldman and C. Moore, Ergodic Equivalence relations,
Cohomology, and Von Neumann Algebras, Trans. of the AMS, V. 234,
P. 289- , 1977.

\bibitem{FL}
W. Feller, An Intorduction to Probabilty Theory and it's
Applications, V. 1, 1968.

\bibitem{FI92}
D. Fiebig and U. R. Fiebig, Covers for Coded Systems, Contemporary
Matematics, Volume 135 P. 139-179, 1992.

\bibitem{HI05}
T. Hamachi and K. Inoue, Embedding of Shifts of Finite Type into
the Dyck Shift, Monatshefte f\"ur Mathematik, Volume 145 P.
107-129, 2005.


\bibitem{PS97}
K. Petersen and K. Schmidt, Symmetric Gibbs Measures, Transactions
of the American Mathematical Society V.349 P. 2775-2811, 1997.


\bibitem{WK74}
W. Krieger, On the Uniqueness of the Equilibruim State,
Mathematical Systems Theory 8, P. 97-104, 1974.


\end{thebibliography}
\end{document}